\documentclass[11pt]{article}
\usepackage{amsfonts}
\usepackage{graphics,amssymb,amsthm,amsmath,epsf}
\usepackage{graphicx}
\usepackage{subfigure}
\usepackage{color}
\usepackage[round]{natbib}
\newtheorem{theorem}{Theorem}[section]
\newtheorem{lemma}{Lemma}[section]

\newtheorem{example}{Example}[section]
\numberwithin{equation}{section}
\newtheorem{remark}{Remark}[section]
\newtheorem{definition}{Definition}[section]

\allowdisplaybreaks
\begin{document}
\title{Cumulative Tsallis Entropy for Maximum Ranked Set Sampling with Unequal Samples}
\author{S. Tahmasebi$^{1}$, M. Longobardi \thanks{maria.longobardi@unina.it} $^{2}$,  M.R. Kazemi$^{3}$, M. Alizadeh$^{1}$ \\
{\small $^{1}$ Department of Statistics, Persian Gulf University , Bushehr, Iran}\\
{\small $^{2}$  Dipartimento di Matematica e Applicazioni Universit$\acute{a}$ di Napoli Federico II Via Cintia, I-80126 Napoli, Italy} \\
{\small $^{3}$ Department of Statistics, Faculty of Science, Fasa University, Fasa, Iran}}
\date{}
\maketitle
\begin{abstract}
In this paper, we consider the information content of maximum ranked set sampling procedure with unequal samples (MRSSU) in terms of Tsallis entropy which is a non-additive generalization of Shannon entropy. We obtain several results of Tsallis entropy including bounds, monotonic properties, stochastic orders, and sharp bounds under some assumptions. We also compare the uncertainty and information content of MRSSU with its counterpart in the simple random sampling (SRS) data. Finally, we develop some characterization results in terms of cumulative Tsallis entropy and {residual} Tsallis entropy of MRSSU and SRS data.


\end{abstract}
\section{Introduction} \label{sec.intro}
The concept of ranked set sampling (RSS) was first introduced by
\cite{mcIntyre-52} to estimate the mean pasture yields and indicated that RSS is a more efficient sampling method in comparison with SRS in terms of the population mean estimation. RSS and some of its variants
have been successfully applied in different areas of applications such as industrial statistics, environmental and ecological studies, biostatistics and statistical genetics. We assume that ${\boldsymbol{X}_{\rm SRS}}=\{X_{i}, i=1,\dots,n\}$ denotes a SRS  of size $n$ from a continuous distribution with probability density function (pdf) $f$ and cumulative distribution function (cdf) $F$.
The one-cycle ranked set sampling involves an initial ranking of $n$ samples of size $n$ as follows:
\[
\begin{array}{ccccccc}
1:       & \underline{\boldsymbol{X_{(1:n)1}}} &  X_{(2:n)1}    & \cdots   & X_{(n:n)1}& \rightarrow & X_{(1)1}=X_{(1:n)1}  \\
2:       & X_{(1:n)2}  & \underline{\boldsymbol{ X_{(2:n)2}}}    & \cdots       & X_{(n:n)2}&\rightarrow & X_{(2)2}=X_{(2:n)2} \\
\vdots   & \vdots    &\vdots     & \ddots       & \vdots &\vdots &\vdots \\
n:       & X_{(1:n)n}&  X_{(2:n)n}        & \cdots  & \underline{\boldsymbol{X_{(n:n)n}}}    &   \rightarrow & X_{(n)n}=X_{(n:n)n}     \\
\end{array}
\]
where $X_{(i:n)j}$  denotes the $i$th order statistic from the $jth$ SRS of size $n$. The resulting sample
is called a  RSS  of size $n$ and denoted by  ${\boldsymbol X}_{\rm RSS}=\{X_{(i)i}, i=1,\dots,n\}$.
   Here,
$X_{(i)i}$ is the $i$th order statistic in a set of size $n$ obtained from the $i$th sample with pdf
\begin{eqnarray*}\label{fXr}
 f_{(i:n)}(x)=\frac{1}{B(i,n-i+1)}f(x)[F(x)]^{i-1}[1-F(x)]^{n-i},
\end{eqnarray*}
where $B(i,n-i+1)$ is the beta function with parameters $i$ and $n-i+1$.

\cite{bir-san-14} proposed  MRSSU  to estimate the mean of the exponential distribution and indicated that MRSSU is better than that of the estimator based on SRS. In the MRSSU, we draw $n$ simple random samples, where the size of the $i$-th samples is $i$, $i=1,...,n$.
 The one-cycle MRSSU  involves an initial ranking of $n$ samples of size $n$ as follows:
\[
\begin{array}{ccccccc}
1:       & \underline{\boldsymbol{X_{(1:1)1}}} &     &    &   & \rightarrow & Y_{1}=X_{(1:1)1}  \\
2:       & X_{(1:2)2}  & \underline{\boldsymbol{ X_{(2:2)2}}}    &       &   &\rightarrow & Y_{2}=X_{(2:2)2} \\
\vdots   & \vdots    &\vdots     & \ddots       & \vdots &\vdots &\vdots \\
m:       & X_{(1:n)n}&  X_{(2:n)n}        & \cdots  & \underline{\boldsymbol{X_{(n:n)n}}}    &   \rightarrow & Y_{n}=X_{(n:n)n}     \\
\end{array}
\]
where $X_{(i:j)j}$ denotes the $i$th order statistic from the $jth$ SRS of size $j$. The resulting sample
is called one-cycle MRSSU of size $n$ and denoted by  ${\boldsymbol Y}_{ MRSSU}=\{Y_{i}, i=1,\dots,n\}$.
Under the assumption of perfect judgment ranking (\cite{Chen-et-al-04}) , $Y_{i}$ has the same distribution
as $X_{(i)i}$ which is the $i$th order statistic (the maximum) in a set of size $i$ obtained from the $i$th sample with probability density function (pdf)
\begin{equation}\label{pdf}
f_{(i)i}(x)=if(x)[F(x)]^{i-1}.
\end{equation}
and distribution function
\begin{equation}\label{cdf}
F_{(i)i}(x)=[F(x)]^{i}.
\end{equation}

In MRSSU, we measure accurately only $n$ maximum order statistics out of $\sum\limits_{i=1}^{n}i=n(n+1)/2$
ranked units and it allows for an increase in set size without introducing too many ranking errors.
\cite{es-di-ta-18} considered information measures of MRSSU in terms of Shannon entropy, R\'enyi entropy and Kullback-Leibler information. \cite{jo-ah-14} explored the notions of uncertainty and information content of RSS data and compared them with their counterparts in SRS data. \cite{th-ja-es-16} obtained some results on residual entropy for ranked set samples. More recently, \cite{ra-qi-19} considered the problem of the information content of RSS data based on extropy measure and the related monotonic properties and stochastic comparisons. However, little works have been done on entropy and we have also not come across any work on entropy properties of RSS or MRSSU  design.\\
Let $X$ denote the lifetime of a system with pdf $f$ and cdf $F$. \cite{shan-48} introduced a measure of uncertainty associated with
$X$ as
\begin{equation}\label{HX}
 H(X)=-\int_{0}^{+\infty} f(x)\log (f(x))dx.
\end{equation}
It is called entropy and  has been used in various branches of statistics and related fields. We refer the reader to \cite{co-th-91} and references  therein for more details. The measure (\ref{HX}) is additive in nature in the sense that for two independent random variables $X$ and $Y$
 \begin{equation}\label{HXY}
 H(X\ast Y)=H(X)+H(Y),
\end{equation}
where $X\ast Y$ denotes the joint random variable. \cite{Tsallis-88} introduced a non-additive generalization of the Shannon entropy which is given by
\begin{equation}\label{eq-SX}
 S_{\alpha}(X)=\frac{1}{1-\alpha}\left[\int_{0}^{+\infty}f^{\alpha}(x)dx-1\right]=\frac{1}{1-\alpha}\left[\int_{0}^{1}f^{\alpha-1}(F^{-1}(u))du-1\right],\; \alpha>0,\; \alpha\neq1,
\end{equation}
where the entropic index $\alpha$ is any real number. Clearly $\lim_{\alpha\rightarrow 1}S_{\alpha}(X)=H(X)$. See \cite{ge-ts-04} for the details of theory and applications of Tsallis entropy.
 Moreover, the Tsallis entropy is a non-additive entropy as for any two independent random variables $X$ and $Y$
\begin{equation}\label{SXY}
 S_{\alpha}(X\ast Y)=S_{\alpha}(X)+S_{\alpha}(Y)+(1-\alpha)S_{\alpha}(X)S_{\alpha}(Y).
\end{equation}
Many applications of Tsallis entropy such as folded proteins, fluxes of cosmic rays, turbulence and many other applications are given in \cite{tsa-bri-04}. Beside the long applications of Tsallis entropy in many applied sciences, as \cite{Wil-Wol-08} stated, there are situations that their uncertainties can be calculated only by Tsallis entropy and the Shannon entropy fails to provide them.The Tsallis entropy has also been extensively used in image processing and signal processing, refer to \cite{to-be-pa-zh-th-02}, \cite{al-es-me-04} and \cite{Yu-Zj-Ho-09}. Considering importance of Tsallis entropy and MRSSU, we try to extend the concept of Tsallis entropy using MRSSU which can be further used in many fields of science. This paper is organized as follows: In Section \ref{sec.tsallis}, we obtain the Tsallis entropy of MRSSU and its comparison with its counterpart under SRS data. Moreover, we provide bounds, monotonic properties, stochastic orders and sharp bound for Tsallis entropy. In Section \ref{sec-cum-res-tsa}, we consider the cumulative Tsallis entropy and residual Tsallis entropy of MRSSU data. Our results include bounds, stochastic ordering and linear transformations. Finally, we conclude the paper in Section \ref{sec-conc}.

\section{Tsallis entropy in RSS and MRSSU data} \label{sec.tsallis}
From (\ref{eq-SX}), the Tsallis  entropy of $\boldsymbol X_{\rm SRS}$ of size $n$ is given by
\begin{eqnarray}\label{SSRS}
 S_{\alpha}(\boldsymbol X_{\rm SRS})&=&\frac{1}{1-\alpha}\left[\int_{0}^{+\infty}...\int_{0}^{+\infty}f^{\alpha}(x_{1})...f^{\alpha}(x_{n})dx_{1}...dx_{n}-1\right]\nonumber\\
 &=&\frac{1}{1-\alpha}\left[\prod_{i=1}^{n}\int_{0}^{\infty}f^{\alpha}(x_{i})dx_{i}-1\right]\nonumber\\
 &=&\frac{1}{1-\alpha}\left([(1-\alpha)S_{\alpha}(X)+1]^{n}-1\right)\nonumber\\
 &=&\frac{1}{1-\alpha}\left(\left[\int_{0}^{1}f^{\alpha-1}(F^{-1}(u))du\right]^{n}-1\right).
\end{eqnarray}
Tsallis entropy associated with the $i$th order statistic from sample of size $n$ is given by \cite{th-ta-ku-15}  as

\begin{eqnarray*}\label{si}
S_{\alpha}(X_{(i:n)})&=& \frac{1}{1-\alpha}\left[\int_{0}^{+\infty}[f_{(i:n)}(x)]^{\alpha}dx-1\right]\\
&=&\frac{1}{1-\alpha}\left[\frac{B(\alpha(i-1)+1,\alpha(n-i)+1)}{(B(i,n-i+1))^{\alpha}}E[f^{\alpha-1}(F^{-1}(W_{i}))]-1\right],\;\alpha>0,\; \alpha\neq1,
\end{eqnarray*}
where $W_{i}$  has the beta distribution with parameters $\alpha(i-1)+1$ and $\alpha(n-i)+1$. Under the perfect ranking assumption, the Tsallis  entropy of $\boldsymbol X_{RSS}$ of size $n$ is given by
\begin{equation}\label{RSS}
 S_{\alpha}(\boldsymbol X_{\rm RSS})=\frac{1}{1-\alpha}\left(\prod_{i=1}^{n}\int_{0}^{\infty}[f_{(i:n)}(x)]^{\alpha}dx-1\right)
 =\frac{1}{1-\alpha}\left(\prod_{i=1}^{n}\left[(1-\alpha)S_{\alpha}(X_{(i:n)})+1\right]-1\right).
\end{equation}

where $X_{(i:n)}$ is the ith order statistic of size $n$ with pdf $f_{(i:n)}(x)$. Similarly, the Tsallis entropy of $\boldsymbol X_{MRSSU}$ of size $n$ is obtained as follows:
\begin{eqnarray}\label{MRSSU}
 S_{\alpha}(\boldsymbol X_{\rm MRSSU})&=&\frac{1}{1-\alpha}\left(\prod_{i=1}^{n}\int_{0}^{\infty}[f_{(i:i)}(x)]^{\alpha}dx-1\right) \nonumber\\
 &=&\frac{1}{1-\alpha}\left(\prod_{i=1}^{n}\left[(1-\alpha)S_{\alpha}(X_{(i:i)})+1\right]-1\right)\nonumber\\
 &=&\frac{1}{1-\alpha}\left(\prod_{i=1}^{n}\left[\int_{0}^{1}i^{\alpha}u^{\alpha(i-1)}f^{\alpha-1}(F^{-1}(u))du\right]-1\right),
\end{eqnarray}
where
\begin{eqnarray}\label{si}
S_{\alpha}(X_{(i:i)})
=\frac{1}{1-\alpha}\left[\frac{i^{\alpha}}{i\alpha-\alpha+1}E[f^{\alpha-1}(F^{-1}(Z_{i}))]-1\right],\;\alpha>0,\; \alpha\neq1,
\end{eqnarray}
and $Z_{i}\sim beta(\alpha(i-1)+1,1)$.
\begin{example}
Let $U$ be a random variable uniformly distributed on $(0,1)$ with $f(F^{-1}(u))=1,\; 0<u<1$. Then we have $$ S_{\alpha}(\boldsymbol U_{\rm MRSSU})=\frac{1}{\alpha-1}\left(1-\frac{(n!)^{\alpha}}{\prod_{i=1}^{n}(1+\alpha(i-1))}\right).$$
\end{example}
\begin{example}
Let $Z$ have exponential distribution with mean $\frac{1}{\theta}$. Thus, $f(F^{-1}(u))=\theta(1-u),\; 0<u<1$, and we have $$ S_{\alpha}(\boldsymbol Z_{\rm MRSSU})=\frac{1}{\alpha-1}\left(1-\theta^{n(\alpha-1)}[(\alpha-1)!]^{n}\prod_{i=1}^{n}\frac{\Gamma(\alpha(i-1)+1)i^{\alpha}}{\Gamma(\alpha i+1)}\right),$$
where $\Gamma(.)$ is the Gamma function.
\end{example}

\begin{theorem}
Let ${\boldsymbol X}_{ MRSSU}$
 be the MRSSU from population $X$ with pdf $f$ and cdf $F$. Then
$[S_{\alpha}(\boldsymbol X_{\rm MRSSU})(1-\alpha)+1]\leq (n!)^{\alpha}[S_{\alpha}(\boldsymbol X_{\rm SRS})(1-\alpha)+1]$   for  $\alpha>0,\;\alpha\neq1$.
\end{theorem}

{\bf Proof.} See Appendix.

Now, we can prove important properties of $S_{\alpha}(\boldsymbol X_{\rm MRSSU})$ using the usual stochastic ordering. For that we present the following definitions:\\
1. The random variable $X$ is said to be smaller than $Y$ according to stochastically ordering (denoted by $X\leq _{st}Y$) if $P(X\geq x)\leq P(Y\geq x)$ for all $x\in \mathbb{R}$, or equivalently $P(X\leq x)\geq P(Y\leq x)$. It is known that  $X \leq _{st}Y \Leftrightarrow \mathbb{E}(\phi(X))\leq \mathbb{E}(\phi(Y))$   for all increasing functions $\phi$ (see \cite{sha-shan-07}).

2. We say that $X$ is smaller than $Y$ in the hazard rate order, denoted by  $X \leq _{hr}Y$, if $\lambda_{X}(x)\geq\lambda_{Y}(x)$ for all $x$.

3. We say that $X$ is smaller than $Y$ in the dispersive order, denoted by  $X \leq _{disp}Y$, if $f(F^{-1}(u))\geq g(G^{-1}(u))$ for all $u\in(0,1)$, where $F^{-1}$ and $G^{-1}$  are right continuous inverses of $F$ and $G$, respectively.

4. A non-negative random variable $X$ is said to have increasing (decreasing) failure rate [IFR (DFR)] if $\lambda(x)=\frac{f(x)}{\bar{F}(x)}$ is increasing (decreasing) in $x$.

5. We say that $X$ is  smaller than $Y$ in the convex transform order, denoted by $X \leq _{c}Y$, if $G^{-1}F(x)$ is a convex function on the support of $X$.

6. A non-negative random variable $X$ is  smaller than $Y$ in the star order , denoted by $X \leq _{*}Y$, if $\frac{G^{-1}F(x)}{x}$ increasing in $x\geq0$.

7. We say that $X$ is  smaller than $Y$ in the supper additive order, denoted by $X \leq _{su}Y$, if $G^{-1}F(t+u)\geq G^{-1}F(t)+G^{-1}F(u)$ for $t\geq 0, u\geq0$.

8. A non-negative random variable $X$ with cdf $F$ is said to have increasing failure rate average (IFRA) if $\frac{\lambda(x)}{x}$ is increasing function in $x>0$. Note that IFR classes of distributions are included to IFRA classes of distributions.

9. A non-negative random variable $X$ with cdf $F$ is new better than used (NBU), if  $\bar{F}(t+u)\leq\bar{F}(t)+\bar{F}(u)$ for $t\geq0$ and $u\geq0$.

\begin{theorem}
Let $X$ and $Y$ be two non-negative random variable with pdf's $f$ and $g$, respectively.  If $X\leq _{disp}Y$, then $S_{\alpha}(\boldsymbol X_{\rm MRSSU})\leq S_{\alpha}(\boldsymbol Y_{\rm MRSSU})$  for $\alpha>0,\;\alpha\neq1$.
\end{theorem}
{\bf Proof.} See Appendix.

\begin{theorem}
If $X\leq _{hr}Y$, and $X$ or $Y$ is DFR, then $S_{\alpha}(\boldsymbol X_{\rm MRSSU})\leq S_{\alpha}(\boldsymbol Y_{\rm MRSSU})$  for $\alpha>0,\;\alpha\neq1$.
\end{theorem}
{\bf Proof.} If $X\leq _{hr}Y$, and $X$ or $Y$ is DFR, then $X\leq _{disp}Y$, due to Bagai and Kochar(1986). Thus, from Theorem (2.2) the desired result follows.

\begin{theorem}
Let $X$ and $Y$ be two non-negative random variable with pdf's $f$ and $g$, respectively. If $X\leq _{su}Y(X\leq _{*}Y \; or \;X\leq _{c}Y)$, then  $S_{\alpha}(\boldsymbol X_{\rm MRSSU})\leq S_{\alpha}(\boldsymbol Y_{\rm MRSSU})$  for $\alpha>0,\;\alpha\neq1$.
\end{theorem}
{\bf Proof.} If $X\leq _{su}Y(X\leq _{*}Y \; or \;X\leq _{c}Y)$, then $X\leq _{disp}Y$, due to Ahmed et al.(1986). So, from Theorem (2.2) the desired result follows.

\begin{theorem}
Let $X$ be a non-negative random variable with decreasing pdf $f$ such that $f(0)\leq1$, then
$$S_{\alpha}(\boldsymbol X_{\rm MRSSU})\geq S_{\alpha}(\boldsymbol U_{\rm MRSSU}),\;\alpha>0,\;\alpha\neq1,$$
\end{theorem}
where $S_{\alpha}(\boldsymbol U_{\rm MRSSU})$ is defined in Example 2.1.\\
{\bf Proof.} The non-negative random variable $X$ has a decreasing pdf if and only if $U\leq_{c} X$, where $U\sim Uniform(0,1)$ (see \cite{sha-shan-07}). Hence, from Theorem (2.4) the desired result follows.

\begin{theorem}
Let $X\in IFR (IFRA,NBU)$, then  $$S_{\alpha}(\boldsymbol X_{\rm MRSSU})\leq S_{\alpha}(\boldsymbol Z_{\rm MRSSU}),\; \alpha>0,\;\alpha\neq1,$$
where $S_{\alpha}(\boldsymbol Z_{\rm MRSSU})$ is defined in Example 2.2.
\end{theorem}
{\bf Proof.} $X\in IFR(IFRZ,NBU)$ if and only if   $X\leq _{c}Z(X\leq _{*}Z \; or \;X\leq _{su}Z)$(see Theorem 4.8.11 of \cite{sha-shan-07}). So, from Theorem (2.4) the desired result follows.

\begin{theorem}
Let $X$ and $Y$ be two independent non-negative random variables. If $X$
and $Y$ have log-concave densities, then
$$S_{\alpha}(\boldsymbol X_{\rm MRSSU}+\boldsymbol Y_{\rm MRSSU})\geq max\{ S_{\alpha}(\boldsymbol X_{\rm MRSSU}),S_{\alpha}(\boldsymbol Y_{\rm MRSSU})\}.$$
\end{theorem}
{\bf Proof.}  See Appendix.

\begin{theorem}
If $f(F^{-1}(u))\geq1,\; 0<u<1$, then $S_{\alpha}(\boldsymbol X^{(n)}_{\rm MRSSU})$ is decreasing(increasing) in $n\geq1$ for $0<\alpha<1(\alpha>1)$.
\end{theorem}
{\bf Proof.} See Appendix.

In the following examples, if we have a system consisting of only two components, then we can compare the Tsallis entropy of $\boldsymbol X_{SRS}$, $\boldsymbol X_{RSS}$ and $\boldsymbol X_{MRSSU}$ of size $n =2$.

\begin{example}\label{exam1}
Let $X$ be uniformly distributed in $(0,b)$, then from (\ref{SXY}), (\ref{si}) and (\ref{SSRS}) for $\boldsymbol X_{SRS}$ , $\boldsymbol  X_{RSS}$ and $\boldsymbol  X_{MRSSU}$ of size $n=2$, respectively, we have
\begin{equation}
\begin{split}
&S_{\alpha}(\boldsymbol X_{SRS})=2 S_{\alpha}(X_{1})+(1-\alpha)(S_{\alpha}(X_{1}))^2=\frac{1}{\alpha-1}[1-b^{2-2\alpha}],\\
&S_{\alpha}(\boldsymbol X_{RSS})=S_{\alpha}(X_{(1:2)})+S_{\alpha}(X_{(2:2)})+(1-\alpha)S_{\alpha}(X_{(1:2)})S_{\alpha}(X_{(2:2)})=\frac{1}{\alpha-1}\left[1-\frac{4^{\alpha}}{(\alpha+1)^2}b^{2-2\alpha}\right],\\
& S_{\alpha}(\boldsymbol X_{MRSSU})=S_{\alpha}(X)+S_{\alpha}(X_{(2:2)})+(1-\alpha)S_{\alpha}(X)S_{\alpha}(X_{(2:2)})=\frac{1}{\alpha-1}\left[1-\frac{2^{\alpha}}{\alpha+1}b^{2-2\alpha}\right].\\
\end{split}
\end{equation}

The differences between $S_{\alpha}(\boldsymbol X_{RSS})$,  $S_{\alpha}(\boldsymbol X_{MRSSU})$  and $S_{\alpha}(\boldsymbol X_{SRS})$ for $n=2$ are given by

\begin{equation*}\label{delta}
\delta_{\alpha}^{(1)}(b)=S_{\alpha}(\boldsymbol X_{RSS})-S_{\alpha}(\boldsymbol X_{SRS})=\frac{b^{2-2\alpha}}{\alpha-1}\left[1-\frac{4^{\alpha}}{(\alpha+1)^2}\right],
\end{equation*}
\begin{equation*}\label{delta}
\delta_{\alpha}^{(2)}(b)=S_{\alpha}(\boldsymbol X_{MRSSU})-S_{\alpha}(\boldsymbol X_{SRS})=\frac{b^{2-2\alpha}}{\alpha-1}\left[1-\frac{2^{\alpha}}{\alpha+1}\right],
\end{equation*}
\begin{equation*}\label{delta}
\delta_{\alpha}^{(3)}(b)=S_{\alpha}(\boldsymbol X_{RSS})-S_{\alpha}(\boldsymbol X_{MSRSU})=\frac{b^{2-2\alpha}}{\alpha-1}\left[\frac{2^{\alpha}}{\alpha+1}-\frac{4^{\alpha}}{(\alpha+1)^2}\right]=
(\frac{2^{\alpha}}{\alpha+1})\delta_{\alpha}^{(2)}(b).
\end{equation*}
In the sequel, Figure \ref{fig-ex-2.1.1} shows the values of $\delta_{\alpha}^{(i)}(b),\; i=1,2,3$ \; for $0< \alpha<1$ and $\alpha>1$. Similarly, in Figure \ref{fig-ex-2.1.2} using $\delta_{\alpha}^{(i)}(b)$ (the difference between Tsallis entropies), we compared the Tsallis entropy of $\boldsymbol X_{SRS}$, $\boldsymbol X_{RSS}$ and $\boldsymbol X_{MRSSU}$ of size $n=3$. Finally, we conclude that $S_{\alpha}(\boldsymbol X_{RSS})\leq S_{\alpha}(\boldsymbol X_{MRSSU})\leq S_{\alpha}(\boldsymbol X_{SRS})$ for $\alpha\neq1$ and $n\geq2$ .
\begin{figure}[ht]
\centering
\includegraphics[width=6.5cm,height=6cm]{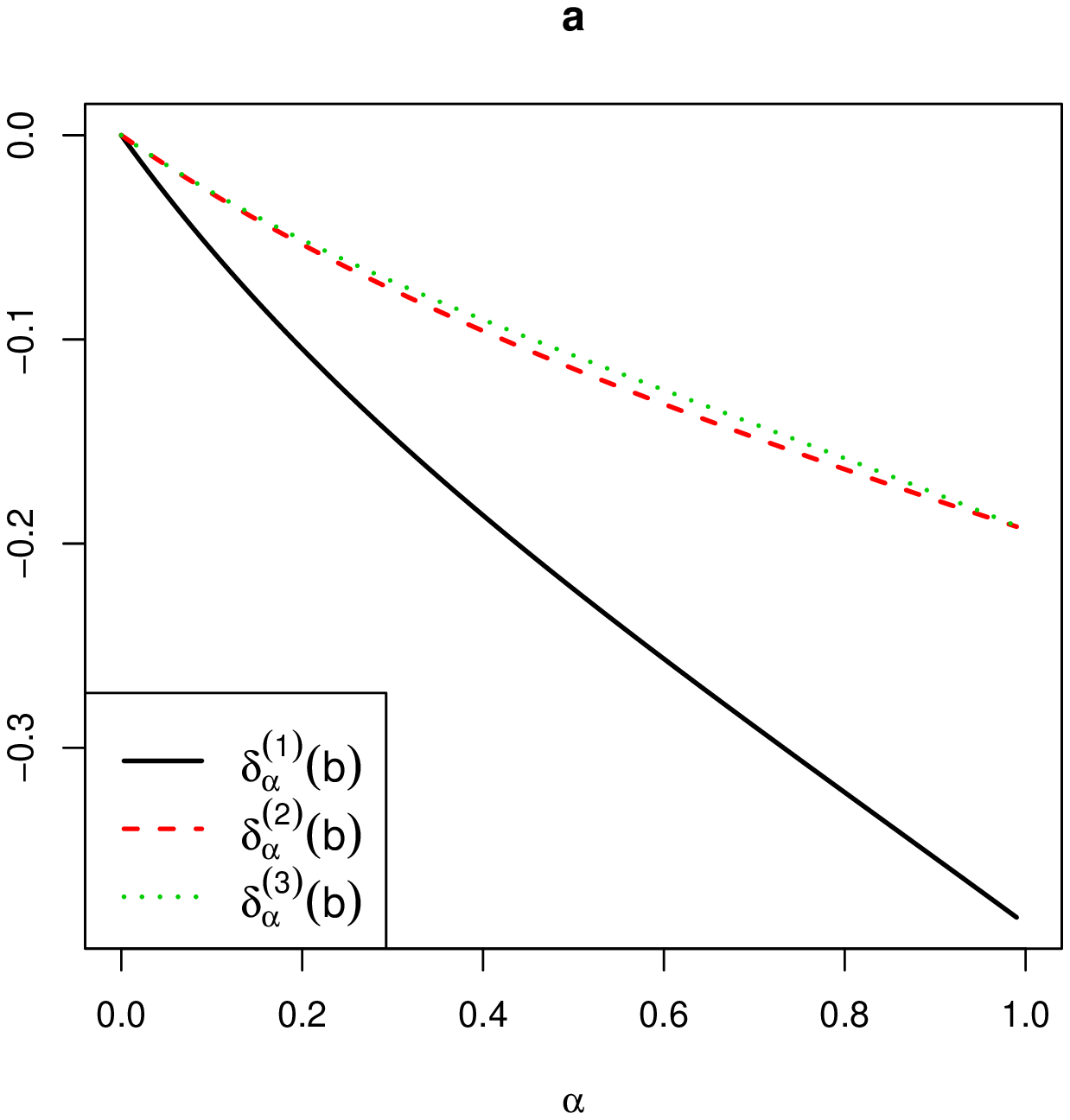}
\includegraphics[width=6.5cm,height=6cm]{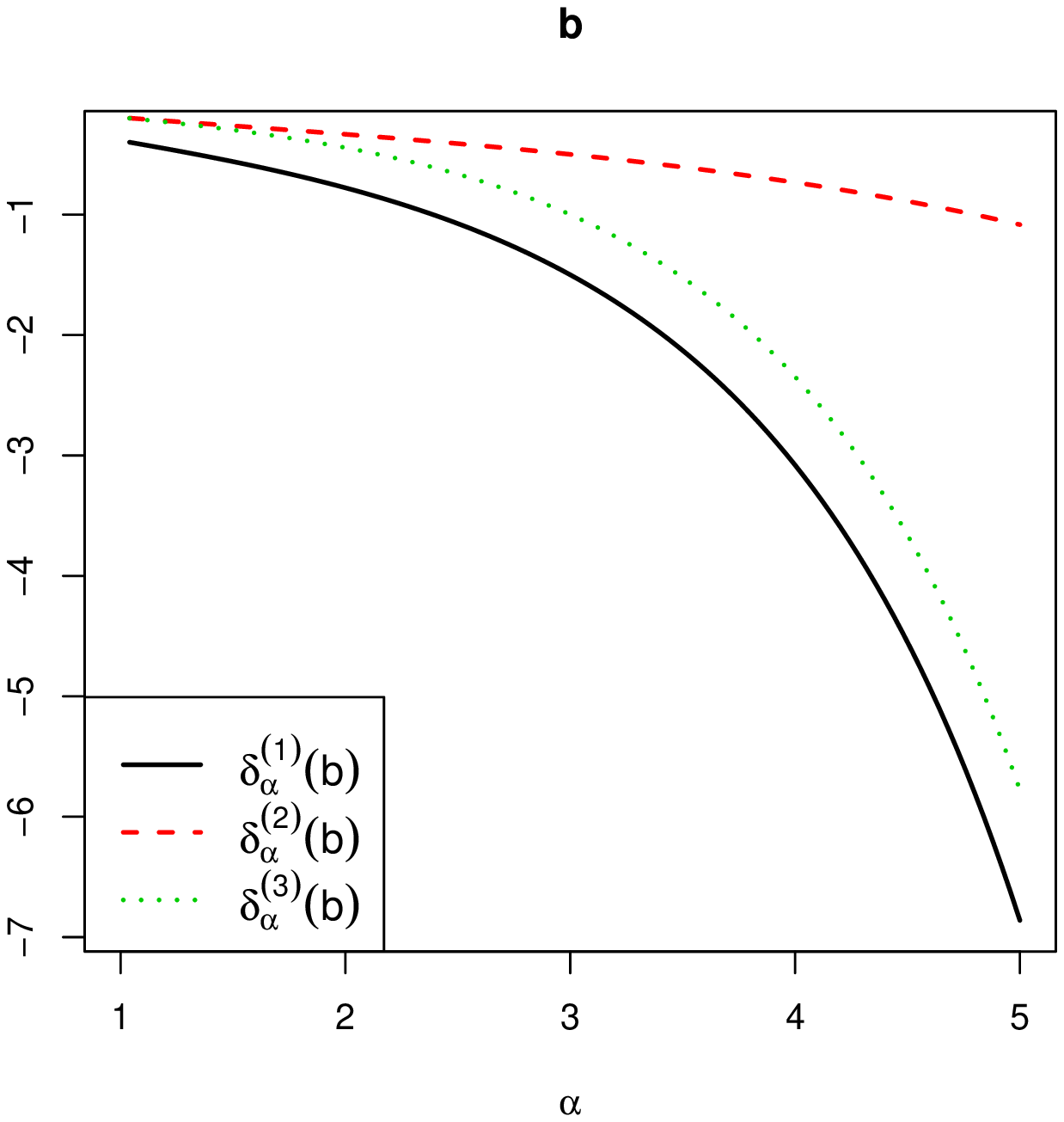}
\vspace{-0.8cm}
\caption{Values of $\delta_{\alpha}^{(i)}(b), i=1,2,3$ for $n=2$ and (a) $0<\alpha<1$ , (b) $\alpha>1$.}
\label{fig-ex-2.1.1}
\end{figure}

\begin{figure}[ht]
\centering
\includegraphics[width=6.5cm,height=6cm]{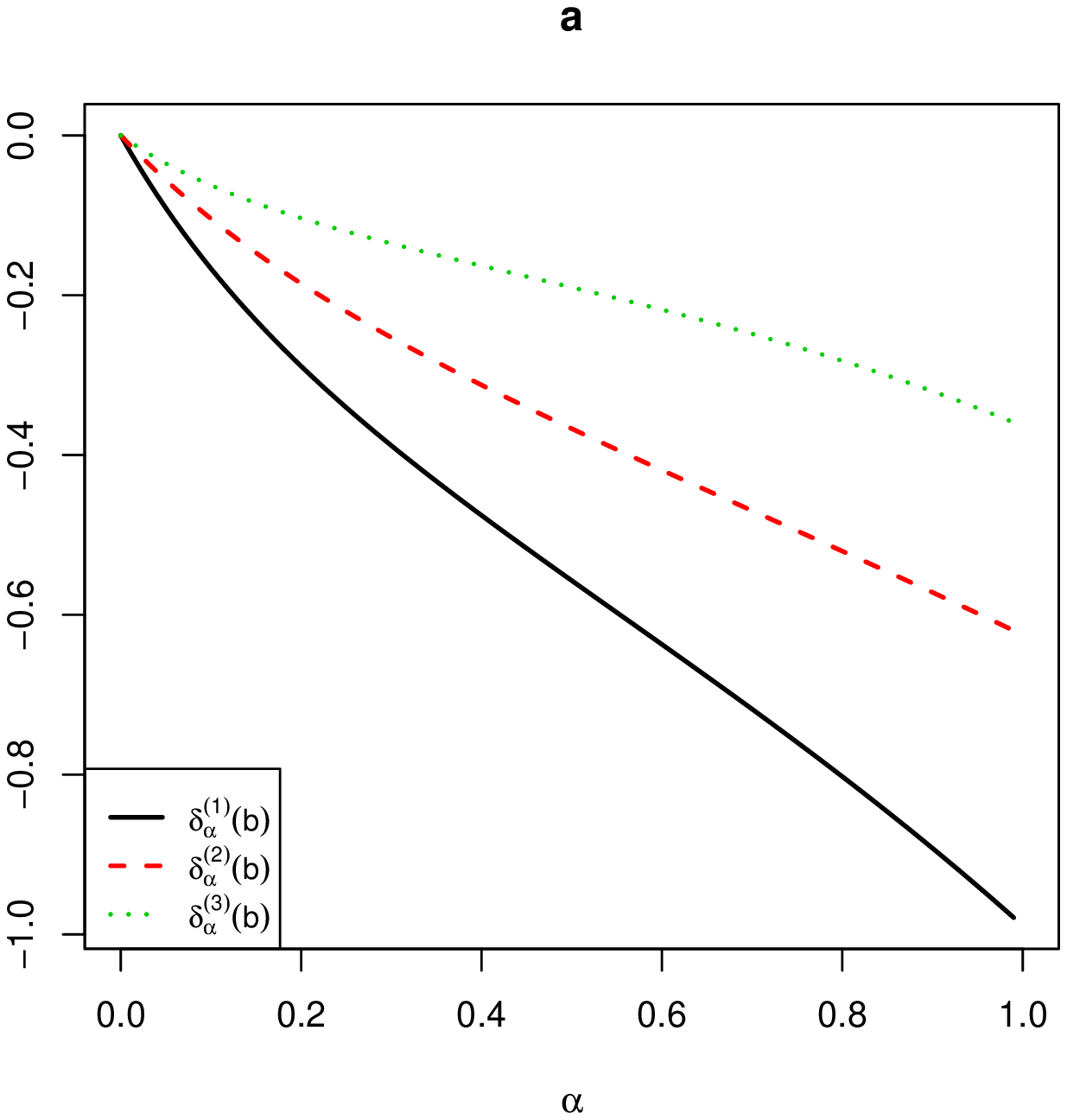}
\includegraphics[width=6.5cm,height=6cm]{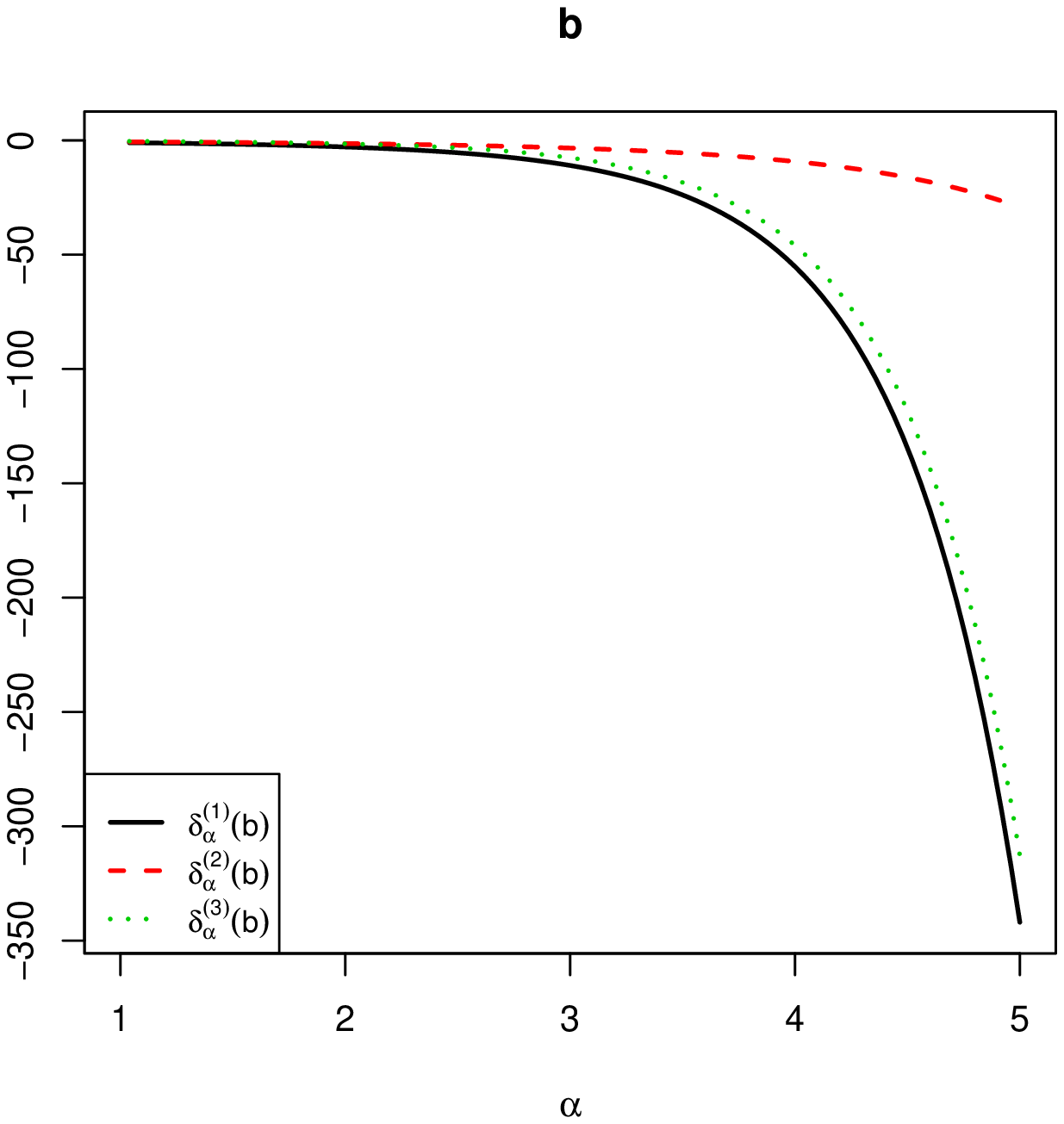}
\vspace{-0.8cm}
\caption{Values of $\delta_{\alpha}^{(i)}(b), i=1,2,3$ for $n=3$ and (a) $0<\alpha<1$ , (b) $\alpha>1$.}
\label{fig-ex-2.1.2}
\end{figure}

\end{example}

\begin{example}\label{Ex2}
Let $X$ be an exponentially distributed random variable with mean $\frac{1}{\theta}$. Straightforward calculations show that for a SRS, perfect RSS and MRSSU of size $n=2$, respectively, we have
\begin{equation*}
\begin{split}
&S_{\alpha}(\boldsymbol X_{SRS})=\frac{1}{\alpha-1}(1-\frac{\theta^{2\alpha-2}}{\alpha^2}),\\
&S_{\alpha}(\boldsymbol X_{RSS})=\frac{1}{\alpha-1}(1-\frac{\theta^{2\alpha-2}4^{\alpha}B(\alpha+1,\alpha)}{2\alpha}),\\
&S_{\alpha}(\boldsymbol X_{MRSSU})=\frac{1}{\alpha-1}(1-\frac{\theta^{2\alpha-2}2^{\alpha}B(\alpha+1,\alpha)}{\alpha}).
\end{split}
\end{equation*}
The difference between  $S_{\alpha}(\boldsymbol X_{RSS})$, $S_{\alpha}(\boldsymbol X_{MRSSU})$  and $S_{\alpha}(\boldsymbol X_{SRS})$ for $n=2$ are given by
\begin{equation}\label{delta2.2.1}
\delta_{\alpha}^{(1)}(\theta)=S_{\alpha}(\boldsymbol X_{RSS})-S_{\alpha}(\boldsymbol X_{SRS})=\frac{\theta^{2-2\alpha}}{(\alpha-1)\alpha^2}\left[1-\frac{4^{\alpha}\alpha B(\alpha+1,\alpha)}{2}\right],
\end{equation}
\begin{equation}\label{delta2.2.2}
\delta_{\alpha}^{(2)}(\theta)=S_{\alpha}(\boldsymbol X_{SRS})-S_{\alpha}(\boldsymbol X_{MRSSU})=\frac{\theta^{2-2\alpha}}{\alpha^{2}(\alpha-1)}\left[2^{\alpha}\alpha B(\alpha+1,\alpha)-1\right],
\end{equation}
\begin{equation}\label{delta2.2.3}
\delta_{\alpha}^{(3)}(\theta)=S_{\alpha}(\boldsymbol X_{RSS})-S_{\alpha}(\boldsymbol X_{MSRSU})=\frac{\theta^{2-2\alpha}2^\alpha B(\alpha+1,\alpha)}{2\alpha(\alpha-1)}\left[2-2^{\alpha}\right].
\end{equation}
From (\ref{delta2.2.1}), (\ref{delta2.2.2}) and (\ref{delta2.2.3}), Figure \ref{fig-ex-2.2.1} shows the values of $\delta_{\alpha}^{(i)}(\theta), i=1,2,3$ for different values of $\alpha$, for simplicity, we take $\theta=1$. Similarly, in Figure \ref{fig-ex-2.2.2} using $\delta_{\alpha}^{(i)}(\theta)$ we compared the Tsallis entropy of $\boldsymbol X_{SRS}$, $\boldsymbol X_{RSS}$ and $\boldsymbol X_{MRSSU}$ of size $n=3$.  It is easily verified that in the exponential distribution case, $S_{\alpha}(\boldsymbol X_{RSS}) \leq S_{\alpha}(\boldsymbol X_{SRS}) \leq S_{\alpha}(\boldsymbol X_{MRSSU})$.

\begin{figure}[ht]
\centering
\includegraphics[width=6.5cm,height=6cm]{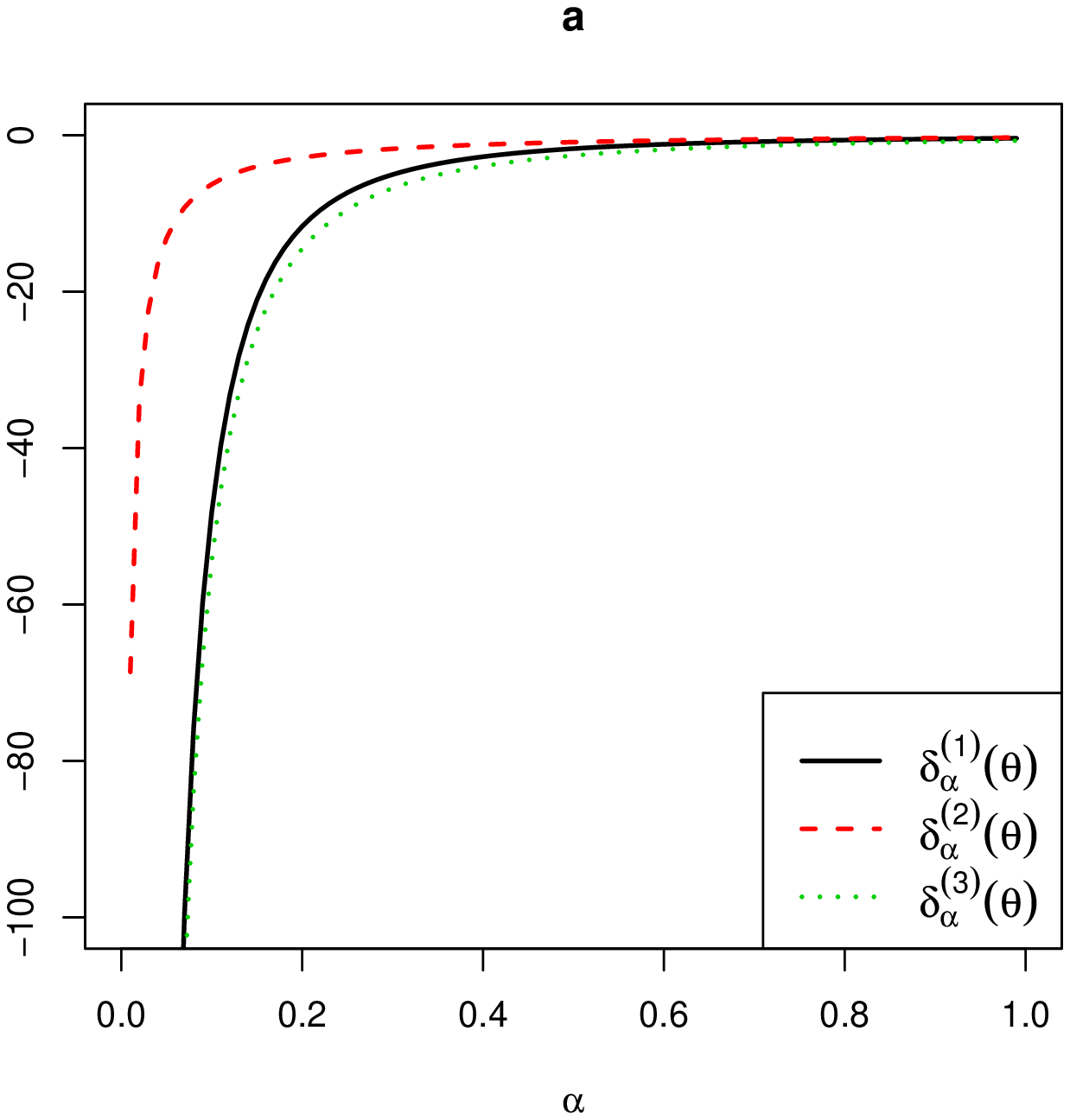}
\includegraphics[width=6.5cm,height=6cm]{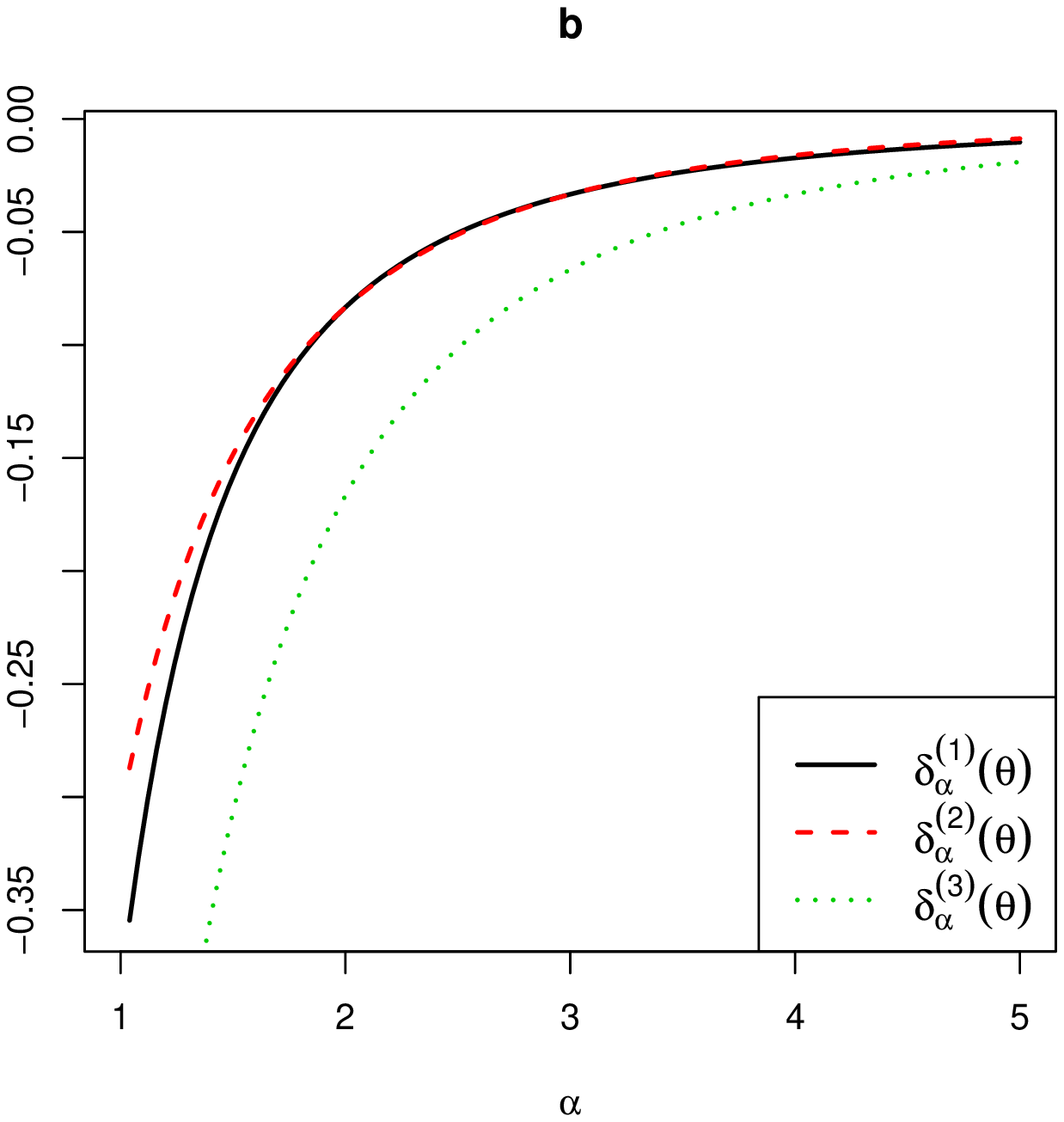}
\vspace{-0.8cm}
\caption{ Values of $\delta_{\alpha}^{(i)}(\theta), i=1,2,3$ for $n=2$ and (a) $0<\alpha<1$ , (b) $\alpha>1$.}
\label{fig-ex-2.2.1}
\end{figure}

\begin{figure}[ht]
\centering
\includegraphics[width=6.5cm,height=6cm]{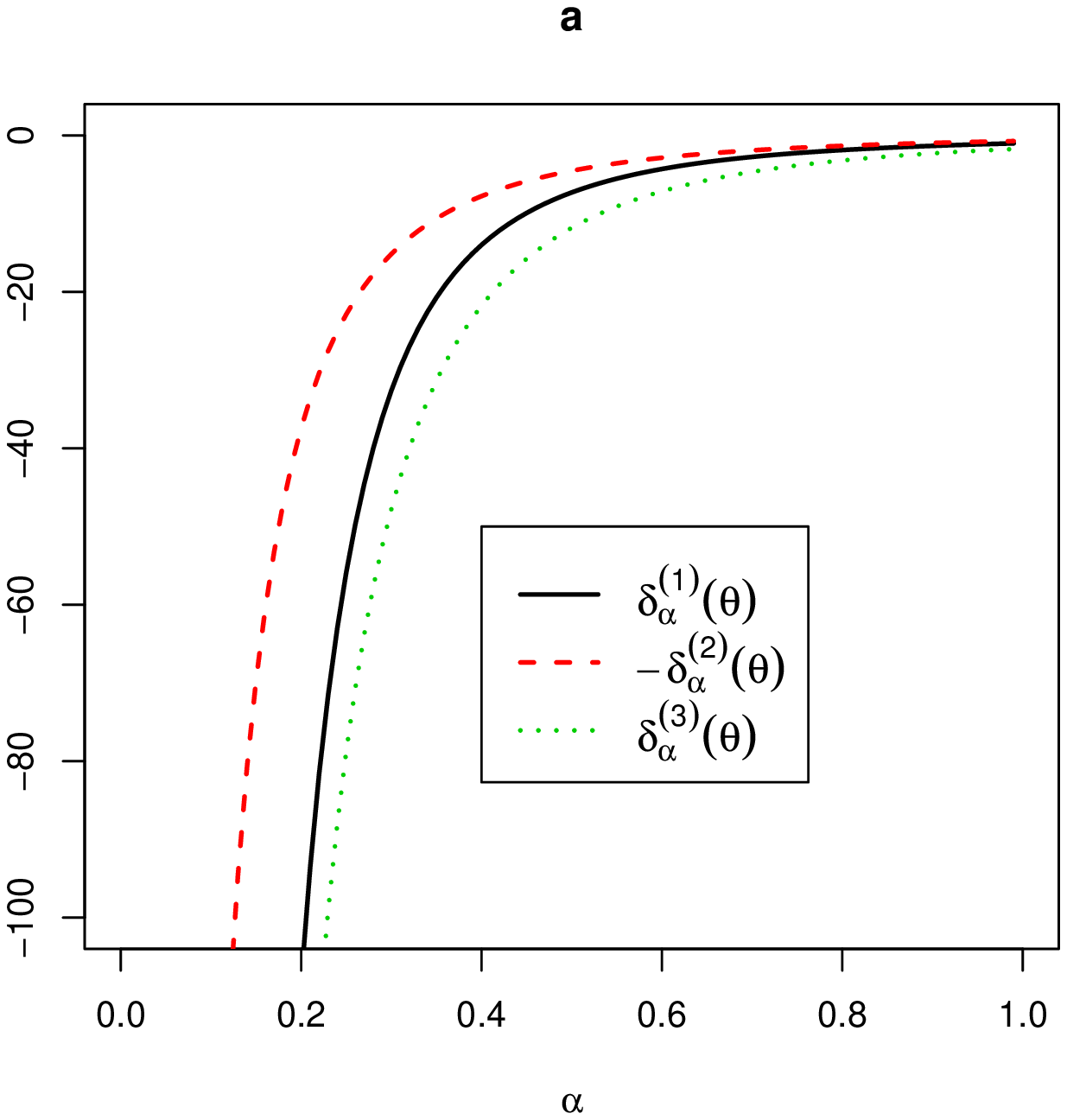}
\includegraphics[width=6.5cm,height=6cm]{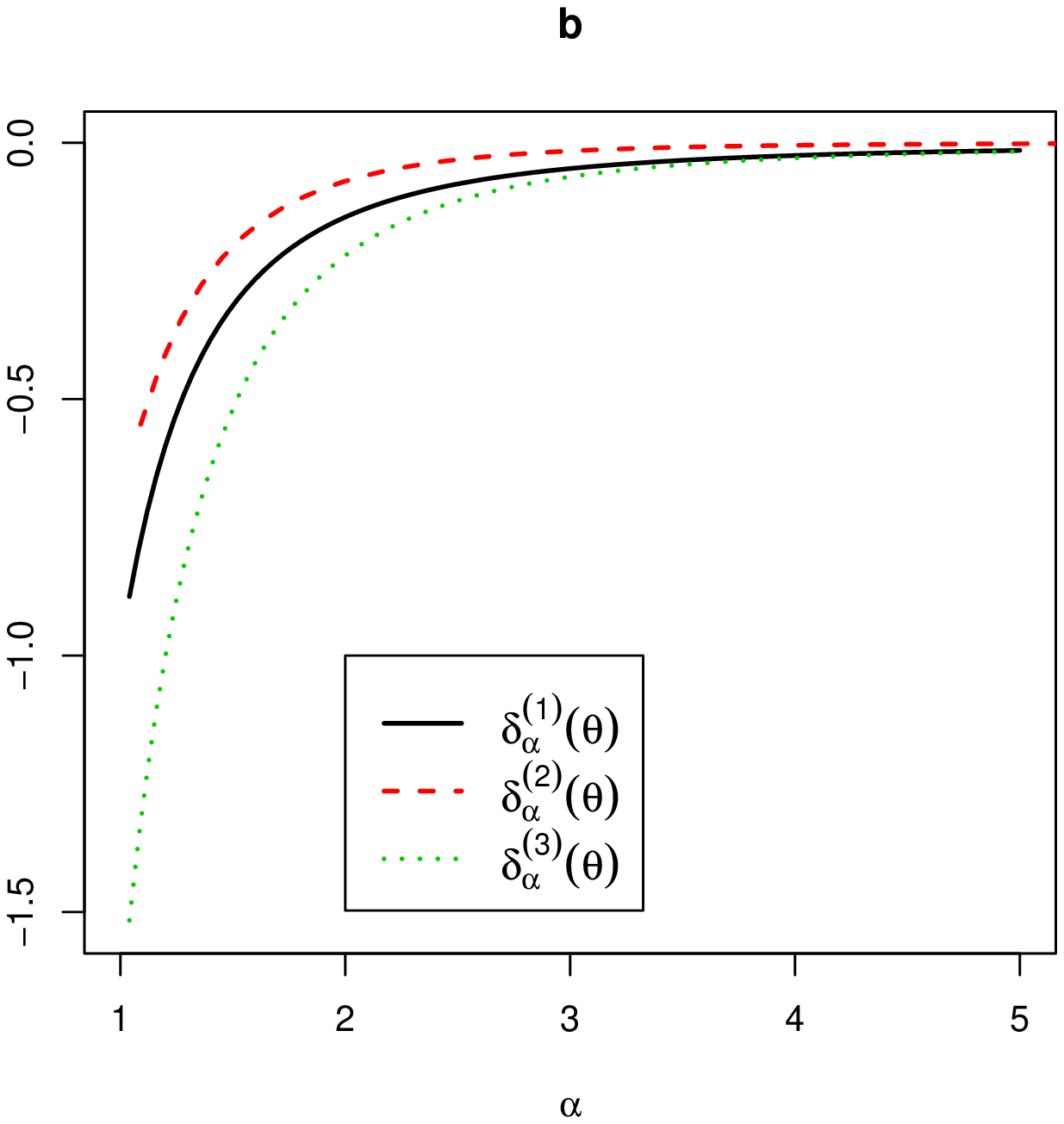}
\vspace{-0.8cm}
\caption{Values of $\delta_{\alpha}^{(i)}(\theta), i=1,2,3$ for $n=3$ and (a) $0<\alpha<1$ , (b) $\alpha>1$.}
\label{fig-ex-2.2.2}
\end{figure}
\end{example}

In the following theorem, we obtain the sharp bounds for $S_{\alpha}(\boldsymbol X_{MRSSU})$ using Steffensen inequalities.
\begin{theorem}
Let $m=\frac{1}{1-\alpha}\left(\prod_{i=1}^{n}\left[\int_{\frac{\alpha(i-1)}{\alpha(i-1)+1}}^{1}i^{\alpha}u^{\alpha(i-1)}f^{\alpha-1}(F^{-1}(u))du\right]-1\right)$,\\
$M=\frac{1}{1-\alpha}\left(\prod_{i=1}^{n}\left[\int_{0}^{\frac{1}{\alpha(i-1)+1}}i^{\alpha}u^{\alpha(i-1)}f^{\alpha-1}(F^{-1}(u))du\right]-1\right)$ and $f$ never increases. Then\\
(i). For  $0<\alpha<1$, we have $m<S_{\alpha}(\boldsymbol X_{MRSSU})<M$.\\
(ii). For  $\alpha>1$, we have $M<S_{\alpha}(\boldsymbol X_{MRSSU})<m$.\\
(iii).  If $f(x)$ never decreases, then all inequalities in parts (i) and (ii) are reversed.
\end{theorem}
{\bf Proof.} See Appendix.

\section{Cumulative Tsallis entropy of MRSSU}\label{sec-cum-res-tsa}
\cite{dic-lon-09} introduced and studied the cumulative entropy (CE) and its dynamic version, which
are defined as
 \begin{eqnarray}
{\mathcal{CE}}(X)=-\int_{0}^{+\infty}F(x)\log F(x)dx,
\end{eqnarray}
and
\begin{eqnarray}
{\mathcal{CE}}(X;t)=-\int_{0}^{t}\frac{F(x)}{F(t)}\log\left( \frac{F(x)}{F(t)}\right)dx,
\end{eqnarray}
respectively. Note that $0\leq{\mathcal{CE}}(X)\leq \infty$. More properties on CE in past lifetime are available in \cite{dic-lon-09}  and \cite{na-del-as-10}. \cite{Ku-18} proposed a cumulative Tsallis entropy (CTE) measure and its dynamic version as
\begin{eqnarray}\label{eq-cte1}
{\mathcal{CE}}_{\alpha}(X)=\frac{1}{1-\alpha}\left(\int_{0}^{+\infty}[F(x)]^{\alpha}dx-1\right),
\end{eqnarray}
\begin{eqnarray}\label{eq-cte2}
{\mathcal{CE}}_{\alpha}(X,t)=\frac{1}{1-\alpha}\left(\int_{0}^{t}\left(\frac{F(x)}{F(t)}\right)^{\alpha}dx-1\right).
\end{eqnarray}

From (\ref{eq-cte1}), the  CTE of $\boldsymbol X_{\rm SRS}$ and $\boldsymbol X_{\rm MRSSU}$ of size $n$ are given by
\begin{equation}\label{CTES}
{\mathcal{CE}}_{\alpha}(\boldsymbol X_{\rm SRS})=\frac{1}{1-\alpha}\left[\left(\int_{0}^{\infty}F^{\alpha}(x)dx\right)^{n}-1\right],
\end{equation}
\begin{equation}\label{CTEMR}
{\mathcal{CE}}_{\alpha}(\boldsymbol X_{\rm MRSSU})=\frac{1}{1-\alpha}\left[\prod_{i=1}^{n}\int_{0}^{\infty}F^{i\alpha}(x)dx-1\right].
\end{equation}
\begin{theorem}
Let ${\boldsymbol X}_{MRSSU}$
be the MRSSU from population $X$ with pdf $f$ and cdf $F$. Then
${\mathcal{CE}}_{\alpha}(\boldsymbol X_{MRSSU})\leq(\geq){\mathcal{CE}}_{\alpha}(\boldsymbol X_{SRS})$   for  $0<\alpha<1(\alpha>1)$.
\end{theorem}
{\bf Proof.} See Appendix.

\begin{theorem}
If $X\leq _{st}Y$, then for $0<\alpha<1(\alpha>1)$ we have $${\mathcal{CE}}_{\alpha}(\boldsymbol X_{MRSSU})\geq(\leq){\mathcal{CE}}_{\alpha}(\boldsymbol Y_{MRSSU}).$$
\end{theorem}
{\bf Proof.} By the assumption on the stochastically ordering, $F^{i\alpha}(x)\geq G^{i\alpha}(x)$ for all $x\geq0$. Now using (3.6) for $0<\alpha<1(\alpha>1)$, we get the desired result.

\begin{theorem}
Let ${\boldsymbol Y}_{ MRSSU}=a{\boldsymbol X}_{ MRSSU}+b$ with $a>0$ and $b\geq 0,$ then we have

$${\mathcal{CE}}_{\alpha}(\boldsymbol Y_{MRSSU})=a^{n}{\mathcal{CE}}_{\alpha}(\boldsymbol X_{MRSSU})+\frac{a^{n}-1}{1-\alpha}.$$
\end{theorem}
{\bf Proof.} The proof is similar to that of Lemma 5.1 of Eskandarzadeh et al. (2018).\\
From (3.5), by independence, the CTE of $\boldsymbol X_{SRS}$ for $n=2$ is given by
\begin{equation}\label{eq-HX-SRS}
 {\mathcal{CE}}_{\alpha}(\boldsymbol X_{SRS})=2{\mathcal{CE}}_{\alpha}(X_{1})+(1-\alpha){\mathcal{CE}}_{\alpha}^{2}(X_{1}).
\end{equation}
In the sequel, due to  (3.6) , under MRSSU data, it is easy to show that for case $n=2$, we have
\begin{equation}\label{eq-HX-MRSSU}
 {\mathcal{CE}}_{\alpha}(\boldsymbol X_{MRSSU})=\sum_{i=1}^{2}{\mathcal{CE}}_{\alpha}(X_{(i)i})+(1-\alpha)\prod_{i=1}^{2}{\mathcal{CE}}_{\alpha}(X_{(i)i}),
 \end{equation}
where $${\mathcal{CE}}_{\alpha}(X_{(i)i})=\frac{1-i\alpha}{1-\alpha}{\mathcal{CE}}_{i\alpha}(X).$$

\begin{example}\label{exam3.1}
Let $X$ be uniformly distributed in $(0,1)$, then from (\ref{eq-HX-SRS}) and (\ref{eq-HX-MRSSU})  for $\boldsymbol X_{SRS}$ and $\boldsymbol  X_{MRSSU}$ of size $n=2$, respectively, we have
\begin{equation}
\begin{split}
&{\mathcal{CE}}_{\alpha}(\boldsymbol X_{SRS})=\frac{\alpha^{2}+2\alpha}{(\alpha^{2}-1)(\alpha+1)},\\
&{\mathcal{CE}}_{\alpha}(\boldsymbol X_{MRSSU})=\frac{2\alpha^{2}+3\alpha}{(\alpha^{2}-1)(2\alpha+1)}.
\end{split}
\end{equation}
The difference between ${\mathcal{CE}}_{\alpha}(\boldsymbol X_{MRSSU})$ and ${\mathcal{CE}}_{\alpha}(\boldsymbol X_{SRS})$ is given by
\begin{equation}\label{delta-ex-3-1}
\delta_{\alpha}={\mathcal{CE}}_{\alpha}(\boldsymbol X_{MRSSU})-{\mathcal{CE}}_{\alpha}(\boldsymbol X_{SRS})=\frac{\alpha}{(\alpha^{2}-1)(2\alpha+1)(\alpha+1)}.
\end{equation}
Formula \eqref{delta-ex-3-1} shows that for $0<\alpha<1$, $\delta_{\alpha}<0$, i.e. Tsallis entropy of ${\mathcal{CE}}_{\alpha}(\boldsymbol X_{MRSSU})$  is smaller than that of ${\mathcal{CE}}_{\alpha}(\boldsymbol X_{SRS})$ and the result will be reversed for $\alpha>1$.
\end{example}

\begin{example}\label{exam3.2}
Suppose $X$ has a beta distribution as beta$(\theta,1)$. From (\ref{eq-HX-SRS}) and (\ref{eq-HX-MRSSU}), respectively, we obtain Tsallis entropy for $\boldsymbol X_{SRS}$ and $\boldsymbol X_{MRSSU}$ of size $n=2$ as follows:
\begin{equation}
\begin{split}
&{\mathcal{CE}}_{\alpha}(\boldsymbol X_{SRS})=\frac{1}{1-\alpha}\left[\frac{1}{(1+\theta\alpha)^2}-1\right],\\
&{\mathcal{CE}}_{\alpha}(\boldsymbol X_{MRSSU})=\frac{1}{1-\alpha}\left[\frac{1}{(1+\theta\alpha)(1+2\theta\alpha)}-1\right].
\end{split}
\end{equation}
The difference between ${\mathcal{CE}}_{\alpha}(\boldsymbol X_{MRSSU})$ and ${\mathcal{CE}}_{\alpha}(\boldsymbol X_{SRS})$ is given by
\begin{equation}\label{delta-ex-3-2}
\delta_{\alpha}={\mathcal{CE}}_{\alpha}(\boldsymbol X_{MRSSU})-{\mathcal{CE}}_{\alpha}(\boldsymbol X_{SRS})=\frac{1}{(\alpha-1)(1+\theta\alpha)^2}\left[1-\frac{1+\theta\alpha}{1+2\theta\alpha}\right].
\end{equation}
It is clear that for $0<\alpha<1$, $\delta_{\alpha}<0$ and for the case $\alpha>1$, $\delta_{\alpha}>0$.
\end{example}

\begin{remark}
The Tsallis entropy in (\ref{eq-SX}) can also be expressed as
\begin{equation}\label{eq-SX-alter}
 S_{\alpha}(X)=\frac{1}{\alpha-1}\left[\int_{0}^{+\infty}(f(x)-f^{\alpha}(x))dx\right].
\end{equation}
Based on (\ref{eq-SX-alter}), recently \cite{cali-et-al-17} introduced an alternate measure of CTE of order $\alpha$ as
\begin{eqnarray}\label{eq-CEX-alter}
{\mathcal{C\xi}}_{\alpha}(X)=\frac{1}{\alpha-1}\left(\int_{0}^{+\infty}(F(x)-(F(x))^{\alpha})dx\right).
\end{eqnarray}
\end{remark}

Due to \eqref{eq-CEX-alter}, under the SRS and MRSSU data , it is easy to show that
\begin{eqnarray}
{\mathcal{C\xi}}_{\alpha}(\boldsymbol X_{SRS})=\frac{1}{\alpha-1}\left[\left(\int_{0}^{+\infty} F(x)dx\right)^{n}-\left(\int_{0}^{+\infty} F^{\alpha}(x)dx\right)^{n}\right],
\end{eqnarray}

\begin{eqnarray}
{\mathcal{C\xi}}_{\alpha}(\boldsymbol X_{MSRSU})=\frac{1}{\alpha-1}\left[\prod_{i=1}^{n}\int_{0}^{+\infty} F^{i}(x)dx-\prod_{i=1}^{n}\int_{0}^{+\infty} F^{i\alpha}(x)dx\right].
\end{eqnarray}

Recently, Abbasnejad and Arghami (2011) defined the following cumulative entropy called failure entropy of order $\alpha$ and its dynamic version as
\begin{eqnarray}\label{eq-cte}
{\mathcal{FE}}_{\alpha}(X)=\frac{1}{1-\alpha}\log \int_{0}^{+\infty}[F(x)]^{\alpha}dx,
\end{eqnarray}
\begin{eqnarray}
{\mathcal{FE}}_{\alpha}(X,t)=\frac{1}{1-\alpha}\log\int_{0}^{t}\left(\frac{F(x)}{F(t)}\right)^{\alpha}dx-1,
\end{eqnarray}
respectively. In the following, we can rewrite the formula of Tsallis cumulative entropy ${\mathcal{CE}}_{\alpha}(X)$  in terms of ${\mathcal{FE}}_{\alpha}(X)$ as
\begin{eqnarray}
{\mathcal{CE}}_{\alpha}(X)=\frac{1}{1-\alpha}\left(e^{(1-\alpha){\mathcal{FE}}_{\alpha}(X)}-1\right).
\end{eqnarray}
Recently, Eskandarzadeh et al. (2018) showed that under the SRS and MRSSU data
$${\mathcal{FE}}_{\alpha}(\boldsymbol X_{SRS})=n{\mathcal{FE}}_{\alpha}(X),\;\;{\mathcal{FE}}_{\alpha}(\boldsymbol X_{MRSSU})=\sum_{i=1}^{n}{\mathcal{FE}}_{\alpha}(X_{i:i})=\sum_{i=1}^{n}\frac{i\alpha-1}{\alpha-1}{\mathcal{FE}}_{i\alpha}(X).$$
Also, they proved that for $0<\alpha<1(\alpha>1)$, ${\mathcal{FE}}_{\alpha}(\boldsymbol X_{MRSSU})\leq(\geq){\mathcal{FE}}_{\alpha}(\boldsymbol X_{SRS})$.
\begin{theorem}
Let ${\boldsymbol X}_{ MRSSU}$
 be the MRSSU from population $X$ with pdf $f$ and cdf $F$. Then
${\mathcal{CE}}_{\alpha}(\boldsymbol X_{MRSSU};t)\leq(\geq){\mathcal{CE}}_{\alpha}(\boldsymbol X_{SRS};t)$   for  $0<\alpha<1(\alpha>1)$.
\end{theorem}
{\bf Proof.} Recalling (\ref{eq-cte}), we have
\begin{eqnarray}
{\mathcal{CE}}_{\alpha}(X;t)=\frac{1}{1-\alpha}\left(e^{(1-\alpha){\mathcal{FE}}_{\alpha}(X;t)}-1\right).
\end{eqnarray}
Therefore, the proof follows of Theorem (5.3) of Eskandarzadeh et al. (2018).\\
In the following theorem, we obtain bounds for ${\mathcal{CE}}_{\alpha}(\boldsymbol X_{MRSSU};t)$ using Hayashi inequalities.

\begin{theorem}
Let $m_{1}=\frac{1}{1-\alpha}\left(A^{n}\prod_{i=1}^{n}\left[\int_{0}^{\lambda}\left(\frac{F(x)}{F(t)}\right)^{\alpha(i-1)}dx\right]-1\right)$,\\
$M_{1}=\frac{1}{1-\alpha}\left(A^{n}\prod_{i=1}^{n}\left[\int_{t-\lambda}^{t}\left(\frac{F(x)}{F(t)}\right)^{\alpha(i-1)}dx\right]-1\right)$ . Then \\
(i). For $0<\alpha<1$, we have $m_{1}<{\mathcal{CE}}_{\alpha}(\boldsymbol X_{MRSSU};t)<M_{1}$.\\
(ii). For $\alpha>1$, we have $M_{1}<{\mathcal{CE}}_{\alpha}(\boldsymbol X_{MRSSU};t)<m_{1}$,\\
where $A=\frac{1}{F^{\alpha}(t)}$ and $\lambda=\frac{1}{A}[(1-\alpha){\mathcal{CE}}_{\alpha}(X;t)+1]$.
\end{theorem}
{\bf Proof.} See Appendix.

Let $X$ and $Y$ be two non-negative random variables with cdfs $F$ and $F^{*}$, respectively. These variables satisfy the proportional reversed hazard rate model(PRHRM) with proportionality constant $\theta>0$ if
\begin{equation}\label{phrh}
F^{*}(x)=[F(x)]^{\theta}, \quad x>0.
\end{equation}
Under the PRHRM \eqref{phrh} we have
\begin{eqnarray}\label{eq-ce-xs}
{\mathcal{CE}}_{\alpha}(\boldsymbol X^{*}_{SRS})=\frac{1-\theta\alpha}{1-\alpha}{\mathcal{CE}}_{\alpha}(\boldsymbol X_{SRS}),
\end{eqnarray}
and
\begin{eqnarray}\label{eq-ce-xs-ii}
{\mathcal{CE}}_{\alpha}(\boldsymbol X^{*}_{MRSSU})=\frac{1-\theta\alpha}{1-\alpha}{\mathcal{CE}}_{\alpha}(\boldsymbol X_{MRSSU}).
\end{eqnarray}

\subsection{Residual Tsallis Entropy}

In the survival analysis and life testing, the current age of the system under consideration is also taken into account. Let $X$ be an absolutely continuous random variable which denotes the lifetime of a system or living organism with pdf $f$. Then $H(X)$ in (\ref{HX}) is not applicable to a system which has survived for some unit of time. \cite{ebrahimi-96} proposed the entropy of the residual lifetime of the random variable $X_{t}=[X-t|X>t]$ as
\begin{equation}\label{eq-ebra-re}
H(X;t)=-\int_{t}^{+\infty}\frac{f(x)}{\bar{F}(t)}\log\left(\frac{f(x)}{\bar{F}(t)}\right)dx, \; t>0.
\end{equation}
 The residual entropy is time-dependent and measures the uncertainty of the residual lifetime of the system when it is still operating at time t. Using \eqref{eq-ebra-re},  the residual Tsallis entropy can be presented by (see, \cite{nan-pau-06} and \cite{Ku-Ta-11})
\begin{equation}\label{sxt}
S_{\alpha}(X;t)=\frac{1}{1-\alpha}\left[\int_{t}^{+\infty}\frac{f^{\alpha}(x)}{\bar{F}^{\alpha}(t)}dx-1\right].
\end{equation}
Also, for any two independent random variables $X$ and $Y$
\begin{equation}\label{sxyt}
S_{\alpha}(X,Y;t)=S_{\alpha}(X;t)+S_{\alpha}(Y;t)+(1-\alpha)S_{\alpha}(X;t)S_{\alpha}(Y;t).
\end{equation}
 In the following, we derive residual Tsallis entropy (\ref{sxt}) for the RSS. Before the main result we obtain residual Tsallis entropy for $i$th order statistics.\\
The residual Tsallis entropy of $i$th order statistics from sample of size $n$ given by
\begin{equation}\label{sxit}
S_{\alpha}(X_{(i:n)};t)=\frac{1}{1-\alpha}\left[\int_{t}^{+\infty}\frac{f^{\alpha}_{(i)}(x)}{\bar{F}^{\alpha}_{(i)}(t)}dx-1\right],
\end{equation}
where, $\bar{F}_{(i)}(t)$ is the survival function of $i$th order statistics and can be represented as
\begin{equation}\label{FI}
\bar{F}_{(i)}(t)=\sum_{j=0}^{i-1}\binom{n}{j}[F(t)]^j [\bar{F}(t)]^{n-j}=\frac{\bar{B}_{F(t)}(i,n-i+1)}{B(i,n-i+1)},
\end{equation}
where $B(a,b)$ and $\bar{B}_{t}(a,b)=\int_{t}^1 u^{a-1}(1-u)^{b-1}du$ are the beta and incomplete beta functions, respectively. For more details about order statistics, one can refer to \cite{ar-ba-na-92}.
Finally by (\ref{sxit}) and (\ref{FI}) we can obtain
\begin{equation}\label{sxitt}
\begin{split}
S_{\alpha}(X_{(i:n)};t)&=\frac{1}{1-\alpha}\left[\frac{1}{\bar{F}_{(i)}(t)}\int_{t}^{+\infty}f^{\alpha}_{(i)}(x)dx-1\right]\\
&=\frac{1}{1-\alpha}\left[\frac{1}{\bar{B}^{\alpha}_{F(t)}(i,n-i+1)}\int_{t}^{+\infty}[f(x)F^{i-1}(x)\bar{F}^{n-i}(x)]^{\alpha}dx-1\right].
\end{split}
\end{equation}
 Under the MRSSU data of size $n=2$, it is easy to show that
\begin{equation}\label{srsst}
S_{\alpha}(\boldsymbol X_{MRSSU};t)=\sum_{i=1}^2 S_{\alpha}(X_{(i:i)};t)+(1-\alpha)\prod_{i=1}^2 S_{\alpha}(X_{(i:i)};t).
\end{equation}
\begin{example}\label{exam4.1}
Suppose that $X$ has a uniform distribution on $(0,1)$. From (\ref{sxyt}) and (\ref{srsst}), respectively, we obtain residual Tsallis entropy for $\boldsymbol X_{SRS}$ and $\boldsymbol X_{MRSSU}$ of size $n=2$,
\begin{equation}
\begin{split}
&S_{\alpha}(\boldsymbol X_{SRS};t)=\frac{1}{1-\alpha}\left[(1-t)^{2-2\alpha}-1\right],\\
&S_{\alpha}(\boldsymbol X_{MRSSU};t)=\frac{1}{1-\alpha}\left[2^{\alpha-1}(1+t)^{1-\alpha}(1-t)^{2-2\alpha}-1\right].
\end{split}
\end{equation}
Let
\begin{equation}\label{dt-ex-4-1}
\tilde{\delta}_{t}=S_{\alpha}(\boldsymbol X_{MRSSU};t)-S_{\alpha}(\boldsymbol X_{SRS};t)=\frac{(1-t)^{2-2\alpha}}{\alpha-1}\left[1-2^{\alpha-1}(1+t)^{1-\alpha}\right].
\end{equation}
Then, from \eqref{dt-ex-4-1}, one can find that $\tilde{\delta}_{t}<0$ for $0<t<1$ when $\alpha>1$, and $\tilde{\delta}_{t}>0$ for $t>1$ when $0<\alpha<1$.
\end{example}
\begin{example}\label{exam4.2}
Suppose that $X$ has an exponential distribution with mean $\frac{1}{\theta}$, by using (\ref{sxyt}) and (\ref{srsst}) for $n=2$ we have
\begin{equation}
\begin{split}
&S_{\alpha}(\boldsymbol X_{SRS};t)=\frac{1}{1-\alpha}\left[\frac{\theta^{2\alpha-2}}{\alpha^2}-1\right]\\
&S_{\alpha}(\boldsymbol X_{MRSSU};t)=\frac{1}{1-\alpha}\left[\frac{\theta^{\alpha-1}2^{\alpha}}{\alpha(2e^{-\theta t}-e^{-2\theta t})^{\alpha}}B_{e^{-\theta t}}(\alpha,\alpha+1)-1\right]
\end{split}
\end{equation}
where $B_{Z}(a,b)=\int_{0}^{z}u^{a-1}(1-u)^{b-1}$ is incomplete beta function. Let $\delta_{t}$ be the difference between $S_{\alpha}(\boldsymbol X_{MRSSU};t)$ and $S_{\alpha}(\boldsymbol X_{SRS};t)$ as follows:
\begin{equation}\label{dt-ex-4-2}
\delta_{t}=S_{\alpha}(\boldsymbol X_{MRSSU};t)-S_{\alpha}(\boldsymbol X_{SRS};t)=\frac{\theta^{2\alpha-2}}{(\alpha-1)\alpha^2}\left[1-\frac{\alpha}{(2e^{-\theta t}-e^{-2\theta t})^{\alpha}}B_{e^{-\theta t}}(\alpha,\alpha+1)\right].
\end{equation}
From (\ref{dt-ex-4-2}), it can be shown that for $0<\alpha<1$, $\delta_{t}<0 \; ; \; \forall t$ and $\delta_{t}>0 \; ; \; \forall t$  when $\alpha>1$.
Following Figure \ref{fig-ex-4-2} shows the behavior of $\delta_{t}$ for different scenarios of parameter $\alpha$.
\begin{figure}[ht]
\centering
\includegraphics[width=6.5cm,height=6cm]{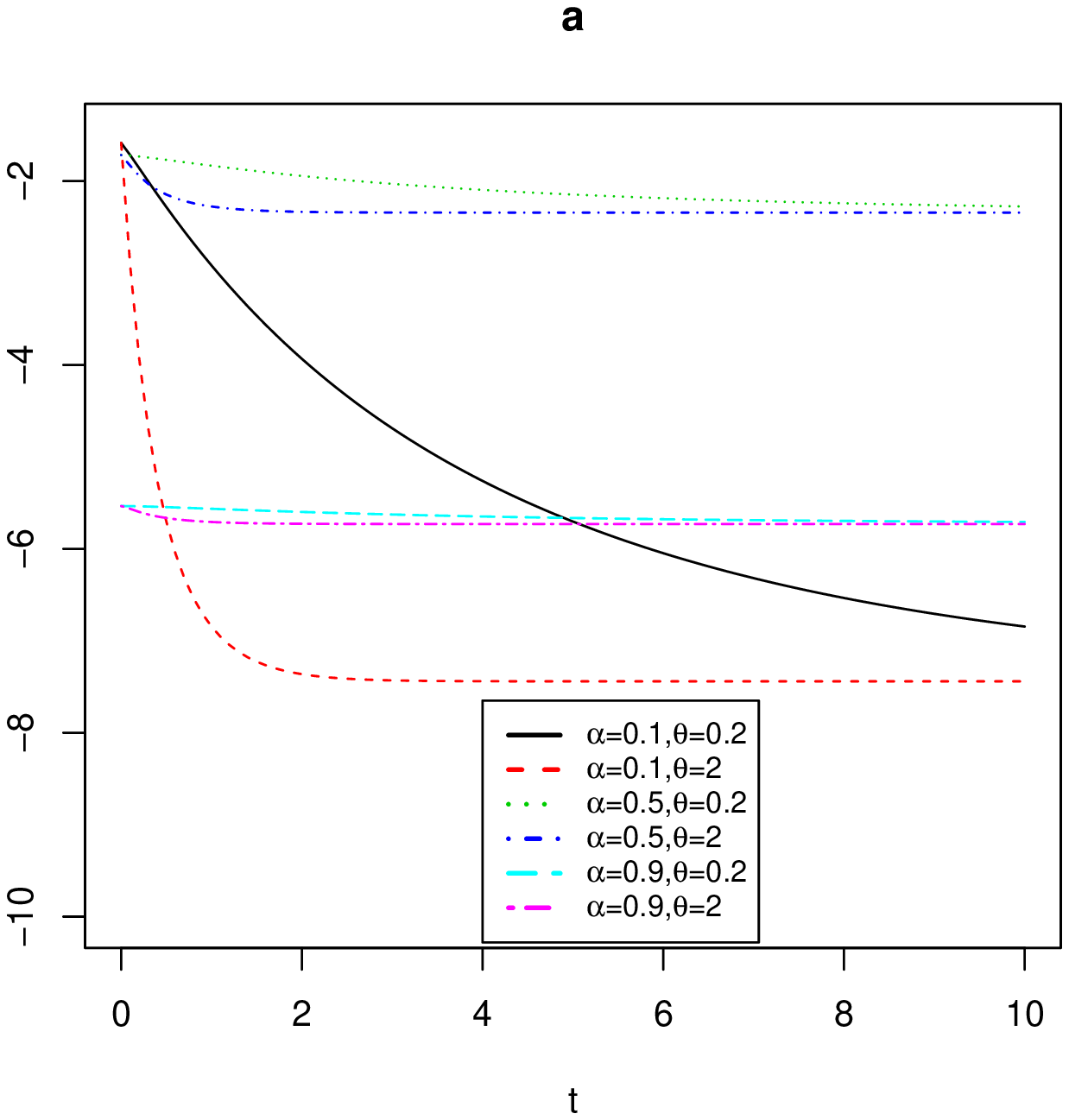}
\includegraphics[width=6.5cm,height=6cm]{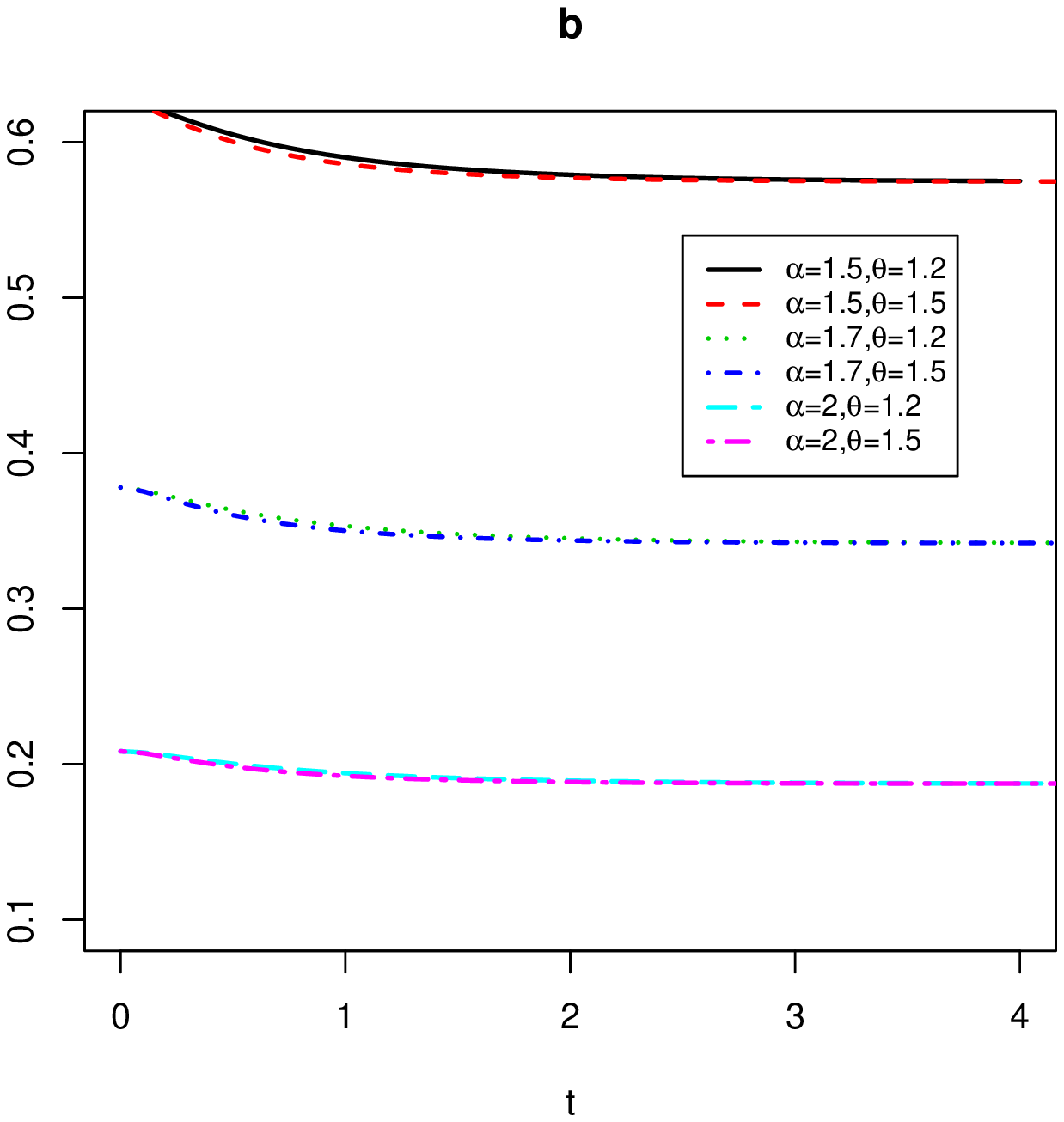}
\vspace{-0.8cm}
\caption{ Values of $\delta_{t}$ for $n=2$ and (a) $0<\alpha<1$ , (b) $\alpha>1$.}
\label{fig-ex-4-2}
\end{figure}
\end{example}
\section{Conclusions} \label{sec-conc}
In this paper, we have considered the information content of MRSSU and SRS data using the Tsallis entropy, cumulative Tsallis entropy and residual Tsallis entropy. We also compared the Tsallis entropy of MRSSU data with SRS and RSS data in the uniform and exponential distributions. For MRSSU data, we obtained several results of Tsallis entropy including bounds, monotonic properties, stochastic orders and sharp bounds under some assumptions. Specifically, we showed that for $0<\alpha<1(\alpha>1)$, ${\mathcal{CE}}_{\alpha}(\boldsymbol X_{MRSSU})\leq(\geq){\mathcal{CE}}_{\alpha}(\boldsymbol X_{SRS})$ and ${\mathcal{CE}}_{\alpha}(\boldsymbol X_{MRSSU};t)\leq(\geq){\mathcal{CE}}_{\alpha}(\boldsymbol X_{SRS};t)$. The results of this paper show some desirable properties of MRSSU compared with the commonly used SRS in the context of the Tsallis entropy and cumulative Tsallis entropy. The concept of ${\mathcal{CE}}_{\alpha}(\boldsymbol X_{MRSSU})$ can be applied in measuring the uncertainty contained in intrinsic fluctuations of  physical systems. Also the Tsallis entropy of $\boldsymbol X_{MRSSU}$ can be used in image or signal processing.
\section{Acknowledgement}
The authors thank the referees and editor for useful comments that improved the paper. Maria Longobardi is supported by GNAMPA (INdAM) and PRIN 2017.

\section{Appendix}
The following definitions and mathematical results will be useful in the computations of the Tsallis entropy for MRSSU.
\begin{definition}
 {\bf Beta function.}
  The Beta function, denoted by $B(a,b)$, is defined as
 \begin{equation}
 B(a,b)=\int_{0}^{1}x^{a-1}(1-x)^{b-1}dx,\;\; a,b>0.
 \end{equation}
\end{definition}
\begin{definition}
 {\bf Gamma function.}
  The Gamma function, denoted by $\Gamma(\alpha)$, is defined as
 \begin{equation}
 \Gamma(\alpha)=\int_{0}^{\infty}x^{\alpha-1}e^{-x}dx,\;\; \alpha>0.
 \end{equation}
\end{definition}

{\bf Proof of Theorem 2.1.} Since
\begin{equation*}
\int_{0}^{1}i^{\alpha}u^{\alpha(i-1)}f^{\alpha-1}(F^{-1}(u))du\leq i^{\alpha}\int_{0}^{1}f^{\alpha-1}(F^{-1}(u))du.
\end{equation*}
The proof follows by recalling \eqref{SSRS} and \eqref{MRSSU}.

{\bf Proof of Theorem 2.2.} By the assumption on the dispersive order, $f(F^{-1}(u))\geq g(G^{-1}(u))$ for all $u\in(0,1)$. Now using \eqref{MRSSU} for $0<\alpha<1(\alpha>1)$, we have
\begin{eqnarray*}
S_{\alpha}(\boldsymbol X_{\rm MRSSU})&=&\frac{1}{1-\alpha}\left(\prod_{i=1}^{n}\left[\int_{0}^{1}i^{\alpha}u^{\alpha(i-1)}f^{\alpha-1}(F^{-1}(u))du\right]-1\right)\nonumber\\
&\leq & \frac{1}{1-\alpha}\left(\prod_{i=1}^{n}\left[\int_{0}^{1}i^{\alpha}u^{\alpha(i-1)}g^{\alpha-1}(G^{-1}(u))du\right]-1\right)=S_{\alpha}(\boldsymbol Y_{\rm MRSSU}).
\end{eqnarray*}

{\bf Proof of Theorem 2.7.}
Suppose that $X$ have a log-concave density, then from Theorem 3.B.7 of \cite{sha-shan-07}, we conclude that $X\leq _{disp} X+Y$ for any random variable $Y$ independent of $X$. Hence, recalling Theorem (2.2),
$ S_{\alpha}(\boldsymbol X_{\rm MRSSU})\leq  S_{\alpha}(\boldsymbol X_{\rm MRSSU}+\boldsymbol Y_{\rm MRSSU})$. Similar result also holds when $Y$ has a log-concave density i.e. $ S_{\alpha}(\boldsymbol Y_{\rm MRSSU})\leq  S_{\alpha}(\boldsymbol X_{\rm MRSSU}+\boldsymbol Y_{\rm MRSSU})$. Therefore, the proof is completed.

\begin{lemma}
Let $\varphi(\alpha)=\frac{(n+1)^{\alpha}}{n\alpha+1}$ be a function of parameter $\alpha$. Then $\varphi(\alpha)$\\
i. has a relative minimum for $0<\alpha<1$.\\
ii. is less than 1 for $0<\alpha<1$.\\
iii. is increasing for  $\alpha>1$.\\
iv. is greater than 1 for $\alpha>1$.
\end{lemma}

{\bf Proof of Theorem 2.8.}
Using \eqref{MRSSU} for $0<\alpha<1(\alpha>1)$, we have
\begin{eqnarray*}
\frac{(1-\alpha)S_{\alpha}(\boldsymbol X^{(n+1)}_{\rm MRSSU})+1}{(1-\alpha)S_{\alpha}(\boldsymbol X^{(n)}_{\rm MRSSU})+1}=(n+1)^{\alpha}\int_{0}^{1}u^{n\alpha}f^{\alpha-1}(F^{-1}(u))du
\leq(\geq)\frac{(n+1)^{\alpha}}{n\alpha+1}.
\end{eqnarray*}
Hence, recalling Lemma (6.1) we have
\begin{eqnarray*}
\frac{S_{\alpha}(\boldsymbol X^{(n+1)}_{\rm MRSSU})}{S_{\alpha}(\boldsymbol X^{(n)}_{\rm MRSSU})}\leq(\geq) 1, \;\; for \;\;0<\alpha<1(\alpha>1).
\end{eqnarray*}
Therefore, the result follows readily.

{\bf Proof of Theorem 2.9 (i).}  Suppose that $f$ never increases. Since $0 \leq u^{\alpha (i-1)} \leq 1$, then we take $\lambda=\int_{0}^{1} u^{\alpha (i-1)}du= \frac{1}{\alpha(i-1)+1}$. Therefore using Steffensen inequalities, we obtain the following inequalities
$$\int_{1-\lambda}^{1}i^{\alpha}f^{\alpha-1}(F^{-1}(u))du\leq \int_{0}^{1}i^{\alpha}u^{\alpha(i-1)}f^{\alpha-1}(F^{-1}(u))du
\leq\int_{0}^{\lambda}i^{\alpha}f^{\alpha-1}(F^{-1}(u))du,
$$
and then we have
\begin{align*}
\prod_{i=1}^{n}\left[\int_{1-\lambda}^{1}i^{\alpha}f^{\alpha-1}(F^{-1}(u))du\right]-1\leq & \prod_{i=1}^{n}\left[\int_{0}^{1}i^{\alpha}u^{\alpha(i-1)}f^{\alpha-1}(F^{-1}(u))du\right]-1\\
 \leq & \prod_{i=1}^{n}\left[\int_{0}^{\lambda}i^{\alpha}f^{\alpha-1}(F^{-1}(u))du\right]-1.
\end{align*}
Hence, the proof is completed. The proof of (ii) and (iii) are similar to proof of part (i).\\

{\bf Proof of Theorem 3.1.}
Since $F^{\alpha}(x)\geq F^{i\alpha}(x)$ for $i\geq1$, we have
\begin{equation*}
\left(\int_{0}^{\infty}F^{\alpha}(x)dx\right)^{n}\geq \prod_{i=1}^{n}\int_{0}^{\infty}F^{i\alpha}(x)dx.
\end{equation*}
The proof follows by recalling (3.5) and (3.6).\\

{\bf Proof of Theorem 3.5(i).} From (\ref{eq-cte2}), we have
$$\lambda=\frac{1}{A}[(1-\alpha){\mathcal{CE}}_{\alpha}(X;t)+1]=\frac{1}{A}\int_{0}^{t}\left(\frac{F(x)}{F(t)}\right)^{\alpha}dx,
 $$
 and
 $$
 \int_{0}^{t}\left(\frac{F(x)}{F(t)}\right)^{i\alpha}dx=\int_{0}^{t}\left(\frac{F(x)}{F(t)}\right)^{\alpha}\left(\frac{F(x)}{F(t)}\right)^{\alpha(i-1)}dx.
 $$
 Since $0\leq\left(\frac{F(x)}{F(t)}\right)^{\alpha}\leq A$ and $\left(\frac{F(x)}{F(t)}\right)^{\alpha(i-1)}$ is never decreases, then using Hayashi inequalities
$$
A \int_{0}^{\lambda}\left(\frac{F(x)}{F(t)}\right)^{\alpha(i-1)}dx\leq
\int_{0}^{t}\left(\frac{F(x)}{F(t)}\right)^{i\alpha}dx
\leq A\int_{t-\lambda}^{t}\left(\frac{F(x)}{F(t)}\right)^{\alpha(i-1)}dx.
$$
Therefore, we have
$$
A ^{n}\prod_{i=1}^{n}\int_{0}^{\lambda}\left(\frac{F(x)}{F(t)}\right)^{\alpha(i-1)}dx\leq
\prod_{i=1}^{n}\int_{0}^{t}\left(\frac{F(x)}{F(t)}\right)^{i\alpha}dx
\leq A^{n}\prod_{i=1}^{n}\int_{t-\lambda}^{t}\left(\frac{F(x)}{F(t)}\right)^{\alpha(i-1)}dx.
$$
Hence, the proof is completed. The proof of (ii) is similar to proof of part (i).\\

\bibliographystyle{apa}

\end{document}